\documentclass[11pt,oneside]{article}

\usepackage{natbib}

\usepackage{graphicx}
\usepackage{amssymb}
\usepackage{amsmath}
\usepackage{amsthm}
\usepackage{indentfirst}

\def\IN{{\mathbb N}}

\def\X{\mathcal{X}}

\def\Var{{\rm Var}}
\def\div{{\rm div}}
\def\grad{\nabla}

\begin{document}
\allowdisplaybreaks
\centerline{\Large\bf
MCMC Confidence Intervals and Biases
}

\bigskip \centerline{by (in alphabetical order)}

\medskip \centerline{Yu Hang Jiang, Tong Liu, Zhiya Lou, Jeffrey S. Rosenthal,
Shanshan Shangguan, Fei Wang, and Zixuan Wu}

\medskip \centerline{\sl Department of Statistical Sciences, University of Toronto}

\medskip\centerline{(December, 2020; last revised June 2021)}

\bigskip
\baselineskip=24pt

\section{Introduction}

Markov chain Monte Carlo (MCMC) is very a powerful tool for estimating
and sampling from complicated high-dimensional distributions (see
e.g.~\citealp{handbook} and the many references therein).  MCMC algorithms
help researchers in a wide spectrum of fields, ranging from Bayesian
statistics to finance to computer science to physics.

One of the biggest challenges when implementing MCMC algorithms
is to evaluate the error of the estimate, which is crucial for
generating accurate results, and can also help when deciding how many
iterations of the chain should be run.  The majority of the existing
results for quantifying MCMC accuracy rely heavily on the Markov
chain Central Limit Theorem (CLT).  However, this CLT is only known
to be valid under specific conditions like geometric ergodicity or
reversibility, which do not always hold and can be difficult to verify
(see e.g.\ \cite{haggstrom2007variance, ibragimov1971independent,
jones2004markov, latuszynski2013nonasymptotic}).  In the reversible
case, \cite{kipnis1986central} established the existence of a CLT for
all reversible Markov chains which have finite asymptotic estimator
variance.  However, without reversibility, CLTs are more challenging.
\cite{Toth1986} generalized the results from \cite{kipnis1986central}
to the non-reversible case, but only under additional conditions which
are very difficult to verify.  And, \cite{Haggstrom2005}
shows that CLTs might not exist for non-reversible chains
under conditions where CLTs would be guaranteed in the reversible case.

There have also been extensive investigations of confidence intervals
in the steady-state literature. These techniques can be split
into three subcategories, outlined in \cite{10.2307/170945}. The most
common approach includes replication and batch means. Replication
runs the simulation numerous times, while batch means divides one
long simulation into batches. Both treat individual runs or batches
as independent and identically distributed random variables. By
assuming each individual group as normally distributed, a confidence
interval can be derived from classical statistical inference
(see details from \citealp{10.1145/256562.256608}). Another approach
is to model the autoregressive structures within the stochastic
process and estimate parameters needed for the confidence interval
from the model. Assuming the process is covariance stationary,
\cite{fishman1971estimating} derives a $p^{\rm th}$ order autoregressive
model as the variance estimator, while an alternative variance
estimator is proposed by \cite{10.1145/358598.358630} using spectral
analysis. Lastly, the regeneration cycles method, first introduced in
\cite{crane1975simulating}, identifies points during the simulation
at which the process ``restarts" probabilistically. Afterwards,
the regeneration epochs are viewed as independent random variables in
order to derive the confidence interval.  However, all of these methods
require assumptions or calculations which can make them difficult to
implement in practice.

The recent paper \cite{rosenthal2017}
derived a simple MCMC confidence interval
which does not require a CLT, using only Chebychev's inequality.
That result required certain assumptions about how the estimator
bias and variance grow with the number of iterations $n$, in particular
that the bias is \(o(1/\sqrt{n})\).
This assumption seemed mild, since
it is generally believed that the estimator bias will be
\(O(1/n)\) and hence \(o(1/\sqrt{n})\) (see e.g.\ page~21 of
\citealp{geyer2011introduction}).
However, questions were raised (\citealp{A1A2tweet}) about how to verify
this assumption, and indeed we show herein (Section~4) that
it might not always hold.


This paper seeks to simplify and weaken the assumptions
in \cite{rosenthal2017}, to make MCMC confidence intervals without
CLTs more widely applicable.  In Section~2,
we derive a conservative asymptotic Markov chain confidence interval
(Theorem~1)
assuming only a finite asymptotic estimator variance as
in \cite{kipnis1986central}, without requiring any bias assumption nor
reversibility nor stationarity nor a CLT;
at significance level $\alpha=0.05$ it
is just 2.3 times as wide as the confidence interval that would follow
from a CLT.  In Section~3, we instead fix the number of iterations $n$,
and obtain corresponding non-asymptotic confidence intervals without
CLT under slightly stronger assumptions.  In Section~4, we consider
the question of when the MCMC bias is or is not \(o(1/\sqrt{n})\),
and show that this property does not always hold but is ensured by a
polynomial ergodicity condition.
In Section~5, we complete the proof of Theorem~1 by extending to
non-stationary chains.
In Section~6, we present some numerical examples to illustrate our results,
and we close in Section~7 with a brief summary.

\section{Asymptotic MCMC Confidence Intervals}


Let \(\{X_n\}\) be a $\phi$-irreducible ergodic Markov chain
on the state space \(\mathcal{X}\),
with stationary distribution $\pi(\cdot)$.
Let \(h: \mathcal{X} \rightarrow
\mathbb{R}\) be a measurable function,
let \(\pi(h) = \int_{x\in\mathcal{X}} \, h(x) \, \pi(x) \, dx\)
be the (finite) expression we wish to estimate, and let
\(e_n := \frac{1}{n}\sum_{i=1}^{n}h(X_i)\) be our estimate.
At significance level $\alpha$, we wish to find
a conservative $1-\alpha$ confidence interval, i.e.\
an interval which contains
$\pi(h)$ with probability at least $1-\alpha$.
Using only a variance bound, we show the following:

\ifx34
(Note this will be
the notation we use throughout the paper.) Define:
\begin{align*}
    \sigma^2 &= \lim_{n\rightarrow\infty}\frac{1}{n}
		E[(\sum_{i=1}^n[h(X_i)-\pi(h)])^2]\\
   &= n\lim_{n\rightarrow\infty}E[(\frac{1}{n}\sum_{i=1}^n[h(X_i)-\pi(h)])^2]\\
    &= n \lim_{n\rightarrow\infty}E[(e_n - \pi(h))^2]\\
    &= n Var(e_n).
\end{align*}
And Bias =$\quad |E(e_n) - \pi(h)|$\\ 
If the chain \(\{X_i\}\) is \(\phi-\)irreducible, reversible, and with
\(\sigma^2 < \infty\), then \(\sqrt{n}-CLT\) holds. We could then use
this fact to derive an asymptotic confidence interval.
\fi

\bigskip\noindent\textbf{Theorem~1:}
If $\limsup_{n \rightarrow \infty} \
n \, Var(e_n) \le B^2$ for  some $B > 0$, then for any $0 < \alpha < 1$ and
$\epsilon > 0$, and $\pi$-a.e.\ initial state $X_0=x\in\X$, the interval
$$
I_{n} := \big(e_n - (1 + \epsilon)n^{-1/2}\alpha^{-1/2}B,
			\ e_n + (1+\epsilon)n^{-1/2}\alpha^{-1/2}B\big)
$$
is an asymptotic conservative $1-\alpha$ confidence interval for $\pi(h)$,
i.e.\
$$
\liminf_{n \to \infty} P[\pi(h) \in I_{n}] \ \geq \ 1 - \alpha.
$$


\begin{proof}
First assume the chain starts in stationarity, so \(E(e_n)=\pi(h)\)
for all \(n \in \mathbb{N}\). Then for any \(a_n > 0\), we
have from Chebychev's inequality that
$$
P\big(|e_n - \pi(h)| \geq a_n\big)
\ = \ P\big(|e_n - E(e_n)| \geq a_n\big)
\ \le \ Var(e_n) \bigm/ a_n^2.
$$
Therefore, setting $a_n = B/\sqrt{n\alpha} > 0$ gives
 \begin{equation*} \begin{split}
 \limsup_{n \to \infty} P\big(|e_n - \pi(h)| \geq a_n\big)
 &\le \limsup_{n \to \infty} \big( Var(e_n) / a_n^2 \big)\\
 &\le \limsup_{n \to \infty} \big(Var(e_n) \frac{n\alpha}{B^2}\big)
 \\ &\le \limsup_{n \to \infty} (B^2 \frac{\alpha}{B^2}) =
 \alpha.
\end{split}\end{equation*}
Then, taking complements gives
\begin{equation*} \begin{split}
 \liminf_{n \to \infty} P\big(
 |e_n - \pi(h)| < n^{-1/2}\alpha^{-1/2}B\big) &=  \liminf_{n
 \to \infty} P\big( |e_n - \pi(h)| < a_n \big) \geq 1-
 \alpha.
\end{split}\end{equation*}
This proves the result (with $\epsilon = 0$)
assuming the chain starts in stationarity.

Finally, applying Theorem~5 from Section~5 below, with $\epsilon>0$
and $r = 1/2$ and $C = \epsilon \, \alpha^{-1/2} B$, we obtain the
result for $\pi$-a.e.\ $X_0=x\in\X$.
\end{proof}

Theorem~1 says that any Markov chain satisfying
\(\limsup_{n\rightarrow\infty} n \, Var(e_n) \leq B^2\) for some \(B >
0\) immediately has a specified asymptotic confidence interval, without
requiring any CLT.  It does not require any bias bound, so it provides a
partial response to \cite{A1A2tweet}.  It does still require a variance
bound.  Asymptotic variance estimators can be obtained in many different
ways, including repeated runs, integrated autocorrelation times, batch
means, window estimators, regenerations, and more (see e.g.\ Section~3
of \citealp{geyer1992}), but they often require challenging conditions
to ensure consistency (\citealp{GlynnWhitt,HJPR,JonesHaran,FlegalJones}).
Alternatively, the variance bound can be consistently
estimated directly from simulations
by sampling $M$ independent copies of the chain and computing the sample
variance of the resulting $e_n$ values
(see the examples in Section~6 below).



For example,
at the usual significance level \(\alpha = 0.05\), taking $\epsilon = 0.001$,
Theorem~1 yields the asymptotic 95\% confidence interval
\[(e_n - 4.48 \, B/\sqrt{n},\ e_n + 4.48 \, B/\sqrt{n}).\]
By contrast, if we knew that a CLT held and that
\(\lim_{n\rightarrow\infty} n Var(e_n) = B^2\),
then we could derive the 95\% confidence interval
\[(e_n - 1.96  \, B /\sqrt{n},\ e_n + 1.96  \, B /\sqrt{n}).\]
The width of the first confidence interval
is 2.3 times the second, but it does not require reversibility,
nor the actual convergence of \(n \, Var(e_n)\) as \(n\to\infty\).

\medskip\noindent\bf Remark. \rm
Our theorems provide {\it conservative} confidence intervals, which
might be larger than necessary, and have coverage
probabilities larger than $1-\alpha$, and lead to running more MCMC
iterations than necessary.  However, we do not consider this to be
a major problem.  The main
challenge of MCMC is to obtain a final answer together with a
guarantee that it is sufficiently accurate.
(For example, \citealp{jones2001honest} seek {\it some} number
of iterations $n'$ for which the chain is
within 0.01 of stationarity, not necessarily the {\it smallest} such $n'$.)
And, conservative confidence intervals do provide such guarantees.
As long as the required
iterations can still be run in a reasonable time, there is no major
disadvantage to running MCMC for somewhat longer than is actually required.

\section{Non-asymptotic MCMC Confidence Intervals}

The confidence intervals from Theorem~1 are only valid asymptotically as
$n\to\infty$.  That limitation is quite common for most MCMC confidence
intervals, since large \(n\) is required for a CLT to hold. However,
since we are not using any CLT in our analysis, it is possible to obtain
a precise non-asymptotic interval, in terms of an upper bound on the
bias, as follows.

\bigskip\noindent
\textbf{Theorem~2:}
Suppose for some fixed $n \in \mathbb{N}$, the chain satisfies the
variance bound $n \, Var(e_n) \le B^2$ for some $B>0$,
and also the bias bound $|E(e_n) - \pi(h)| \le C$ for some $C \ge 0$.
Then for any significance level $\alpha \in (0,1)$,
setting \(\delta = \frac{C}{\frac{B}{\sqrt{n\alpha}} + C} \in [0,1)\)
and \(a_n = \frac{B}{\sqrt{n\alpha}(1-\delta)}\),
the fixed-$n$ interval $$I_{n} \ := \ (e_n - a_n, \ e_n + a_n)$$
is a non-asymptotic conservative $1-\alpha$ confidence interval, i.e.\
$$P[\pi(h) \in I_{n}] \ \geq \ 1 - \alpha, \qquad n\in\IN.$$


\begin{proof} We first compute that
\[\frac{\delta}{1-\delta}
\ = \ \frac{C}{\frac{B}{\sqrt{n\alpha}}+C}
\bigg/ \frac{\frac{B}{\sqrt{n\alpha}}}{\frac{B}{\sqrt{n\alpha}}+C}
\ = \ \frac{C\sqrt{n\alpha}}{B},\]
and hence
$$
\delta a_n
\ = \ \frac{\delta}{1 - \delta} \ \frac{B}{\sqrt{n}\alpha}
\ = \ \frac{C\sqrt{n\alpha}}{B} \ \frac{B}{\sqrt{n\alpha}}
\ = \ C.
$$
Thus $|E(e_n - \pi(h))| \le C = \delta a_n$, and hence
\(a_n - |E(e_n)-\pi(h)| \geq a_n - C = (1-\delta) a_n > 0\).
Therefore, using the triangle inequality and then Chebyshev's inequality,
we have that
\begin{align*} P\bigg(|e_n - \pi(h)| \geq a_n|\bigg) &= P\bigg(|e_n - E(e_n) + E(e_n) - \pi(h)| \geq a_n\bigg) \\ &\le P\bigg(|e_n - E(e_n)| + |E(e_n) - \pi(h)| \geq a_n\bigg) \\&= P\bigg(|e_n - E(e_n)| \geq a_n - |E(e_n) - \pi(h)|\bigg)
\\ &\leq Var(e_n) \bigm/ \big(a_n - |E(e_n) - \pi(h)| \big)^2
\\ &\leq Var(e_n) \bigm/ \big( (1 - \delta) a_n \big)^2
\\ &= Var(e_n) \big(\frac{n \alpha}{B^2}\big)
	\leq B^2\frac{\alpha}{B^2 } = \alpha.
\end{align*}
Taking complements gives $P[\pi(h) \in I_{n}] \geq 1 - \alpha$,
as claimed.
\end{proof}

If the chain is in stationarity, or at least reaches stationarity within
$n$ iterations, then the bias is zero, and we obtain:

\bigskip
\noindent\textbf{Corollary \(\mathbf{1}\):}  Let \(n
\in \mathbb{N}\) be a fixed time
such that the chain is in stationarity after \(n\)
steps. Then if $n \, Var(e_n) \le B^2$ for some $B > 0$,
then for any significance level $0 < \alpha < 1$,
the interval $$I_{n} := (e_n - n^{-1/2}\alpha^{-1/2}B, \ e_n +
n^{-1/2}\alpha^{-1/2}B)$$
is a non-asymptotic conservative $1-\alpha$ confidence interval, i.e.\
$$P[\pi(h) \in I_{n}] \ \geq \ 1 - \alpha.$$

\begin{proof} 
By the stationarity assumption, $|E(e_n) - \pi(h)| = 0$, so we can
apply Theorem~2 with \(C = 0\).  It follows
that $\delta = 0$ and $a_n = n^{-\frac{1}{2}}\alpha^{-1/2}B$.
The result then follows immediately from Theorem~2.
\end{proof}

The assumptions for the non-asymptotic bound in Theorem~2 and Corollary~1
are stronger than for the asymptotic bound of Theorem~1, since they
require a bound on the bias or for the chain to be at stationarity
after \(n\) iterations. However, we will see in the next section that
we can sometimes utilize properties such as polynomial ergodicity to
help us establish a bound on bias. Also, in practice, MCMC users often
approximately verify stationarity through a plethora of convergence
diagnostics such as plots and renewal theory and non-parametric tests;
see e.g.\ \cite{mengersen1999mcmc} for a review.

Next, we present a result which does not assume stationarity, nor
require a bound on the bias, nor require a bound on the variance. But as
a trade-off, it assumes a bound on an absolute first moment, which might
be harder to verify. It could still be useful if e.g.\ the first moment
condition can be linked to another property that is easier to satisfy,
which could be explored in future research.

\bigskip \noindent\textbf{Theorem~3:}
Suppose for some fixed $n \in \mathbb{N}$
we have $E(|e_n - \pi(h)|) \le \gamma_n$ for some constant $\gamma_n>0$.
Then for any significance level $\alpha\in(0,1)$,
the interval $$I_{n} := (e_n -
\gamma_n\alpha^{-1}, \ e_n + \gamma_n\alpha^{-1})$$
is a non-asymptotic conservative $1-\alpha$ confidence interval, i.e.\
$$P[\pi(h) \in I_{n}] \ \geq \ 1 - \alpha.$$

\medskip
\begin{proof} 
Setting $a_n = \gamma_n/\alpha > 0$, we have
by Markov's inequality that
\begin{equation*}\begin{split}
P(|e_n - \pi(h)| \geq a_n)
&\le E(|e_n - \pi(h)|) \bigm/ a_n \\
&= E(|e_n - \pi(h)|) \, \frac{\alpha}{\gamma_n} \\
&\le ( \gamma_n\frac{\alpha}{\gamma_n}) = \alpha.
\end{split}\end{equation*}
Taking complements,
\begin{equation*} \begin{split}
P[\pi(h) \in I_{n}] &= P\big( |e_n - \pi(h)| \le \gamma_n \alpha^{-1} \big) \\
&=  P\big( |e_n - \pi(h)| \le a_n \big) \geq 1- \alpha.
\qedhere
\end{split}\end{equation*}
\end{proof}

In particular, if the \(\gamma_n\) converge monotonically to zero,
then we obtain a confidence interval which shrinks to a point
as \(n\) approaches infinity.

\section{The Order of MCMC Bias}

Since we are estimating the quantity \(\pi(h) =
\int_{x\in\mathcal{X}} h(x) \, \pi(x) \, dx\) by the Markov chain estimator
\(e_n := \frac{1}{n}\sum_{i=1}^{n}h(X_i)\),
the bias after $n$ iterations is given by
$\textbf{Bias}(e_n) := {E}(e_n) - \pi(f)$.
As previously mentioned, the results in
\cite{rosenthal2017} assumed this bias was
\(o(1/\sqrt{n})\) since it is generally believed to be $O(1/n)$
(see e.g. p.~21 of
\citealp{geyer2011introduction}).  However, this is not always the case:

\bigskip
\noindent\textbf{Example~1:}
Consider the Markov chain with state space
\(\mathcal{X} = \{0, 1, 2, 3, \dots\}\), and transition probabilities
given by \(p_{0,0} = 1\),
and for all $n \ge 1$, $p_{n, 0} = 1 - \frac{\sqrt{n}}{\sqrt{n +
1}}\) and \(p_{n, n+1} =\frac{\sqrt{n}}{\sqrt{n + 1}}\). Then the
chain is \(\phi\)-irreducible (and aperiodic) with \(\pi(x) = \phi(x)
= \delta_0(x)\), i.e.\ \(\pi(0) = 1\) and \(\pi(x) = 0\) for all
\(x \neq 0\).
Assume $X_0 = 1$.
We then compute that, for $n=1,2,3,\ldots$,
$$
P[X_n \neq 0]
= P[X_n = n+1]
= \prod_{i=1}^{n}  \frac{\sqrt{i}}{\sqrt{i + 1}} = \frac{1}{\sqrt{n + 1}}.
$$
Thus,
$$
\lim_{n\rightarrow\infty}P[X_n \neq 0]
= \lim_{n\rightarrow \infty}\frac{1}{\sqrt{n + 1}} = 0.
$$
So, the chain will converge to \(\pi(\cdot)\) (from any initial distribution).

Next, consider the function on $\X$ defined by \(f(0) = 0\),
and $f(x) = 1$ for $x \ge 1$.  Then $\pi(f) = f(0) = 0$.
It follows that
\begin{align*}
    \textbf{Bias}(e_n) &= {E}(e_n) - \pi(f) =  {E}(e_n)
	= \frac{1}{n} \sum_{j = 1}^n {E}[h(X_j)] \\
    &= \frac{1}{n} \sum_{j = 1}^n \Big[ f(j + 1)P(X_{j}= j + 1)
					+f(0)P(X_j=0) \Big] \\
    &= \frac{1}{n} \sum_{j = 1}^n \frac{1}{\sqrt{j + 1}}
    	\geq \frac{1}{n} \sum_{j = 1}^n
			\frac{1}{\sqrt{n + 1}}= \frac{1}{\sqrt{n + 1}}.
\end{align*}
On the other hand,
\begin{align*}
   \textbf{Bias}(e_n) &=\frac{1}{n} \sum_{j = 1}^n \frac{1}{\sqrt{j + 1}}
    	\le \frac{1}{n} \int_0^n x^{-\frac{1}{2}} \, dx \\
    &= \frac{1}{n} \, 2\, x^{\frac{1}{2}} \Big|_{x=0}^{x=n}
    	= \frac{1}{n} (2 \sqrt{n}) = \frac{2}{\sqrt{n}}.
\end{align*}
That is, ${1 \over \sqrt{n+1}} \le \textbf{Bias}(e_n) \le {2 \over
\sqrt{n}}$.  In particular, the bias is \(O(1/\sqrt{n})\), but is not
$O(1/n)$ nor even \(o(1/\sqrt{n})\).
\qed
\bigskip

Example~1 raises the question of what conditions guarantee the bias
to be \(o(1/\sqrt{n})\).  We shall derive such a result for a class
of Markov chains that are polynomially ergodic, defined as follows
(see e.g.\ \cite{jarner2003}, \cite{jones2004markov}):

\bigskip
\noindent\textbf{Definition:} Let \(\{X_n\}\) be a Markov chain with
stationary distribution $\pi(\cdot)$,
and let \(||\cdot||\) be total variation distance. Then the chain is {\it
polynomially ergodic} if there exists a function \(M:\X\to[0,\infty)\) such that:
\[||P^n(x, \cdot)-\pi(\cdot)|| \ \leq \ M(x) \, n^{-m},
\quad x\in\X, \ x\in\IN;\]
here \(m > 0\) is the {\it order} of the polynomial ergodicity.

\bigskip \noindent\textbf{Theorem~4:}
Let \(\{X_n\}\) be a polynomially
ergodic Markov chain of order $m > \frac{1}{2}$, with stationary
distribution $\pi(\cdot)$.
Suppose for some \(D \in [0, \infty)\) and \(f: \mathcal{X}
\rightarrow \mathbb{R}\), we have \(|f(x)|\leq D\).
Then for any fixed initial state \(X_0 = x\),
the absolute bias $|\mathbf{Bias}(e_n)|$ is $o(n^{-1/2})$ as $n\to\infty$.

\begin{proof}
Let $X_0=x$.  We compute that
\begin{align*}
        |\mathbf{Bias}(e_n)| &= |{E}(e_n) - \pi(f)|\\
    &\leq \frac{1}{n}\sum_{i=1}^n|{E}[f(X_i)] - \pi(f)|\\
    &\leq \frac{1}{n}\sum_{i=1}^n \sup_{g:\mathcal{X}\rightarrow\mathbf{R}, \ |g(x)|\leq D}|{E}[g(X_i)] - \pi(g)| \\
    &\leq \frac{1}{n}\sum_{i=1}^n 2D \, ||P^i(x, \cdot)-\pi(\cdot)||
\\
    &\leq \frac{1}{n}\sum_{i=1}^n 2DM(x) \, i^{-m} \\
    &= \frac{2DM(x)}{n}\sum_{i=1}^n i^{-m}.
    \end{align*} 
Case 1: $\frac{1}{2} < m < 1$. Then
    \begin{align*} 
     |\mathbf{Bias}(e_n)| &\leq \frac{2DM(x)}{n}\sum_{i=1}^n i^{-m}\\& \leq \frac{2DM(x)}{n}\int_0^n x^{-m}dx
    \\&= \frac{2DM(x)}{n}\cdot \frac{1}{1-m}(n^{1-m}-0^{1-m})\\
    &= \frac{2DM(x)}{n}\cdot \frac{1}{1-m}(n^{1-m})\end{align*} 
Therefore
     $$n^{1/2} |\mathbf{Bias}(e_n)| \leq \frac{2DM(x)}{1-m}n^{1/2-m},$$
which $\to 0$ as $n\to\infty$ since
$m>1/2$ and $1-m>0$ and $D, \, M(x) < \infty$.

Case 2: $m \geq 1$. Find some $\beta$ such that $1/2 < \beta < 1 \le m$. Then since $\sum_{i = 1}^n i^{-m} \le \sum_{i = 1}^n i^{-\beta}$, it follows as
above that
\begin{equation*}
\lim_{n\rightarrow\infty}n^{1/2} |\mathbf{Bias}(e_n)|\le \lim_{n\rightarrow\infty}\frac{2DM(x)}{1-\beta}n^{1/2-\beta}  = 0.
\qedhere
\end{equation*}
\end{proof}

\bigskip\noindent
\textbf{Remark:} This result says the bias is $o(1/\sqrt{n})$
for any polynomially ergodic chain of order more than $1/2$.
In the context of Theorem~2, this means that
we can always find a constant \(C\) such that $|E(e_n) - \pi(h)|
\leq C = \delta a_n$, since $a_n$ is $O(1/\sqrt{n})$. Furthermore,
if we can calculate an explicit value for
\(M(x)\), then we can obtain a value for \(C\).
As a specific example, if a chain has polynomial order $3/4=:m$,
with initial state $X_0 =: x$ satisfying $M(x)=2$,
and variance bound $n \, Var(e_n) \le 4 =: B^2$,
and functional bound $|f(x)| \le 5 =: D$, then after $n=100$ iterations
we will have
$|\mathbf{Bias}(e_n)| \le {2DM(x) \over n} \, {1 \over 1-m} \, n^{1-m}
= (20/n) (4) (n^{1/4}) \doteq 2.53 =: C$, so we can apply Theorem~2
at significance level $\alpha=0.05$ to find that
$\delta = {2.53 \over (2/\sqrt{5}) + 2.53}
\doteq 0.74$ and $a_n = 2/(\sqrt{5} \, (1-0.74)) \doteq 3.44$, giving
the 95\% confidence interval $(e_n-3.44, \, e_n+3.44)$ after 100 iterations.

\section{Extending to Non-Stationary Chains}

Theorem~1 above was initially proved assuming the chain started in
stationarity.  However, in practice MCMC is hardly ever started in
stationarity, so accuracy bounds without this assumption are much
more useful.  We now prove a general result which says that asymptotic
confidence intervals from stationarity can always be enlarged slightly
to become asymptotic confidence intervals from arbitrary initial states.

\bigskip
\noindent\textbf{Theorem~5:} Consider an ergodic Markov chain $\{X_n\}$
on a state space $\X$ with stationary distribution $\pi(\cdot)$,
functional $h$, and usual estimator $e_n$.  Suppose the sequence
$(e_n+a_n,e_n+b_n)$ is an asymptotic conservative $1-\alpha$ confidence
interval for $\pi(h)$ when started in stationarity, i.e.\
\[\liminf_{n \to \infty} P\big(a_n < \pi(h) - e_n < b_n\big)
\ \geq \ 1 - \alpha, \quad X_0 \sim \pi(\cdot).\]
Then for any $c > 0$ and $0 < r < 1$,
and $\pi$-a.e.\ initial state $x\in\X$,
the sequence $(e_n + a_n - cn^{-r}, \, e_n + b_n + cn^{-r})$
is an asymptotic conservative $1-\alpha$ confidence interval
for the chain when started from the initial state $X_0=x$, i.e.\
\[\liminf_{n \to \infty}
P\big(a_n - cn^{-r} < \pi(h) - e_n < b_n + cn^{-r} \big)
	\ \geq \ 1 - \alpha, \quad X_0=x.\]

\begin{proof}
By ergodicity, for $\pi$-a.e.\ $x\in\X$, we have
$\lim_{n\to\infty} ||P^n(x, \cdot) - \pi(\cdot)|| = 0$.  Hence,
for fixed $\epsilon>0$ and $x\in\X$, we can
find $m \in \mathbb{N}$ such that $||P^m(x, \cdot) - \pi(\cdot)||
\le \epsilon$.
Let $\{X_n\}$ be our original chain with $X_0=x$, and let
$\{X'_n\}$ be a second copy of the chain in stationarity, i.e.\
with $X'_0 \sim \pi(\cdot)$ and hence
$X'_n \sim \pi(\cdot)$ for all $n$.
By Proposition 3(g) of \cite{roberts2004general},
we can couple $\{X_ n\}$ and $\{X_n'\}$ such that
$P(H) \ge 1-\epsilon$, where
$$
H \ = \ \{X_n = X_n' \text{ for all } n  \geq m\} .
$$

Now, let
$e_n =\frac{1}{n}\sum_{i = 1}^n h(X_i) $, and
$e_n' = \frac{1}{n}\sum_{i = 1}^n h(X'_i)$ be the estimators
from the two chains, so by assumption we have
\[\liminf_{n \to \infty} P\big(a_n < \pi(h) - e'_n < b_n\big)
\ \geq \ 1 - \alpha.\]
Then on the event $H$, for any $n \geq m$ we have
$$
\Big|\big( \pi(h) - e_n \big) - \big( \pi(h) - e'_n \big)\Big|
\ = \ \frac{1}{n}\Big|\sum_{i = 1}^m \big(h(X'_i) - h(X_i)\big)\Big|
\ =: \ \frac{1}{n} |Z|,
$$
where $Z = \big|\sum_{i = 1}^m \big(h(X'_i) - h(X_i)\big)\big|$.
Hence, if $H$ holds and $a_n < \pi(h) - e'_n < b_n$
and $\frac{1}{n} |Z| \le cn^{-r}$, then
$a_n - cn^{-r} < \pi(h) - e_n < b_n + cn^{-r}$.
Furthermore, $Z$ is a fixed
finite random variable, so there is $A<\infty$ with $P(|Z|>A) \le \epsilon$.
It follows that for $n \ge (A/c)^{1/(1-r)}$, we have
$$
P\big(\frac{1}{n} |Z| > c n^{-r}\big)
\ = \ P\big(|Z| > c n^{1-r}\big)
\ \le \ P\big(|Z| > A\big)
\ \le \ \epsilon
.
$$
We conclude that for all $n \ge \max[m, \ (A/c)^{1/(1-r)}]$,
\begin{equation*}
\begin{split}
&
P\big( \{a_n - cn^{-r} < \pi(h) - e_n < b_n + cn^{-r}\}^C \big)
\\
&\le \
P\big( \{a_n < \pi(h) - e'_n < b_n\}^C \big)
	\, + \, P(H^C) \, + \ P(|Z|>A) \\
&\le \
P\big( \{a_n < \pi(h) - e'_n < b_n\}^C \big)
	\, + \, \epsilon \, + \epsilon,
\end{split}
\end{equation*}
i.e.\
\begin{equation*}
P\big( a_n - cn^{-r} < \pi(h) - e_n < b_n + cn^{-r} \big)
\ \ge \
P\big( a_n < \pi(h) - e'_n < b_n \big)
	- 2 \, \epsilon.
\end{equation*}
Then, taking $\liminf$ gives
\begin{equation*}
\liminf_{n\to\infty}
P\big( a_n - cn^{-r} < \pi(h) - e_n < b_n + cn^{-r} \big)
\ \ge \ \alpha - 2 \, \epsilon.
\end{equation*}
Since this is true for any $\epsilon>0$, we must actually have
\begin{equation*}
\liminf_{n\to\infty}
P\big( a_n - cn^{-r} < \pi(h) - e_n < b_n + cn^{-r} \big)
\ \ge \ \alpha,
\end{equation*}
giving the result.
\end{proof}

\section{Numerical Examples}

In this section, we apply Theorem~1 to various non-reversible
examples, to obtain confidence intervals directly without the need to
establish a CLT nor any convergence rates.

\subsection{A Cyclical Non-Reversible Chain}

Define a Markov chain on the state space \(\mathcal{X} =
\mathbb{N}\), as follows. Fix $0<r<1$.
For all \(x > 1\), let \(p_{1, x} = r^{x-2} \, (1-r)\),
and \(p_{x, x-1} = 1\), with $p_{x,y}=0$ otherwise.
It is easily verified that this
chain has stationary distribution given by:
\begin{align*}
    \pi(1) &= \pi(2) = \frac{1 - r}{r^2 -2r + 2} \, , \\
    \pi(x) &= \Big(\frac{1 - r}{r^2 -2r + 2}\Big) \, r^{x-1} , \quad x > 2.
\end{align*}
Furthermore, this chain is irreducible and aperiodic.
Hence, the chain will converge asymptotically to $\pi$,
and furthermore $e_n$ will converge to $\pi(h)$ whenever $\pi|h|<\infty$.
However, it is not trivial to establish a confidence interval
for $\pi(h)$ in terms of $e_n$ by means of a CLT, since
this chain is clearly non-reversible, and establishing a condition like
geometric ergodicity would require careful drift function arguments.
Instead, we use our method.
Take $X_0=1$, \(r = 0.75\), and \(h(x) = x^{0.5}\).
For different numbers
of iterations \(n\), we run \(M = 100\) replications
to estimate $n \, \Var(e_n)$ and $e_n$.  Our results are as follows:

\begin{center}
\begin{tabular}{ |c|c|c|c| } 
\hline
\(n\) & \(n \, $\Big.$ \widehat{Var}(e_n)\) & \(\widehat{e_n}\) \\
\hline
1,000 & 2.840 & 1.91171 \\ 
2,000 & 3.413 & 1.90454 \\ 
5,000 & 4.054 & 1.90518 \\ 
10,000 & 2.557 & 1.90171 \\ 
20,000 & 3.299 & 1.90384 \\ 
50,000 & 2.611 & 1.90501 \\ 
100,000 & 3.514 & 1.90434 \\ 
200,000 & 4.323 & 1.90414 \\ 
500,000 & 3.040 & 1.90433 \\ 
1,000,000 & 3.596 & 1.90416 \\ 
\hline
\end{tabular}
\end{center}

\noindent
These replications indicate that
\(n \, \Var(e_n) \le 5\) for all $n$.
Hence, we can take \(B = \sqrt{5}\). Then, setting
\(\alpha = 0.05\) and \(\epsilon = 0.001\), Theorem~1 gives
an approximate conservative 95\% confidence interval
for $\pi(h)$ equal to:
\[
I \ \equiv \ I_{1,000,000} \ = \ (1.894, 1.914)
\, .
\]
This is quite a small interval, of width $0.02$,
which provides good confidence about $\pi(h)$,
without needing to establish any CLT or difficult ergodicity property.

\ifx34
\section{Proof for Stationary Distribution (Appendix / To be Omitted)}
First realize that \(\sum_{x=2}^\infty f(x) = (1-r)\sum_{x = 2}^\infty r^{x-2} = 1\). Therefore,
\begin{align*}
    \pi(1) &= \pi(2)\\
    \pi(x) &= \pi(1)f(x) + \pi(x + 1) \quad , x>2.
\end{align*}
By solving this recurrence relation, we get \[    \pi(x) = \pi(1)(1 - \sum_{x=2}^\infty f(x)) = \pi(1)r^{x-1} \quad , x>2.\]
Afterwards, normalize the distribution to get
\begin{align*}
    1 &= \sum_{x = 1}^\infty \pi(x) = \pi(1)(2 + \sum_{x=3}^\infty r^{x-1})\\
    1 &= \pi(1)\frac{r^2 - 2r + 2}{1 - r} \\
        \pi(1) &= \pi(2) = \frac{1 - r}{r^2 -2r + 2}\\
    \pi(x) &= (\frac{1 - r}{r^2 -2r + 2})r^{x-1} \quad , x > 2.
\end{align*}
\fi

\subsection{A Diffusive Non-Reversible Chain}

Consider the Markov chain with state space $\mathcal{X} = \{0, 1,
2, 3, \dots\}$, and transition probability given by $p_{0, 0} =
0.99$, $p_{0, 1} = 0.01$, and for all $x \geq 1$, $p_{x, x + 1} =
\left(\frac{x}{x + 1}\right)^2$ and $p_{x, 0} = 1 - \left(\frac{x}{x
+ 1}\right)^2$.   This chain is easily computed to have stationary distribution:
\begin{equation}\pi(0) = \frac{1}{1 + \frac{0.01 \pi^2}{6}}, \quad
{\rm and} \quad \pi(x)
= \frac{0.01}{x^2} \, \pi(0) \ \text{for} \ x \geq 1.\end{equation}
This chain is again irreducible and aperiodic, so again
$e_n$ will converge to $\pi(h)$ whenever $\pi|h|<\infty$.
However, this chain is again clearly non-reversible, and it is again
non-trivial to establish a CLT to obtain a confidence interval for
$\pi(h)$ in terms of $e_n$.
Instead, we again use our method.
Let $$h(x) = \begin{cases} 0, & x = 0 \\ x^{-1}, & x \geq 1 \end{cases}$$
and again take $X_0=1$.
For different numbers of iterations $n$, we run $M = 1,000$
replications to estimate $e_n$ and its variance.  The results are as follows:
\begin{center}
\begin{tabular}{ |c|c|c|c| } 
\hline
\(n\) & \(n \, $\Big.$ \widehat{Var}(e_n)\) & \(\widehat{e_n}\) \\
\hline
100    & 0.0181   & 0.01400\\ 
200    & 0.0163   & 0.01264 \\ 
500    & 0.0148   & 0.01224 \\
1,000   & 0.0156   & 0.01196 \\ 
2,000   & 0.0151   & 0.01188 \\ 
5,000   & 0.0151   & 0.01195 \\ 
10,000  & 0.0142   & 0.01185 \\ 
20,000  & 0.0153   & 0.01184\\ 
50,000  &  0.0154   & 0.01186 \\ 
100,000 & 0.0141   &  0.01182   \\  
\hline
\end{tabular}
\end{center}
These simulations indicate that \(n \, \Var(e_n) \le 0.02\)
for all $n$, so we can take
\(B = \sqrt{0.02}\). Then, setting
\(\alpha = 0.05\) and \(\epsilon = 0.001\), Theorem~1 gives
an approximate conservative 95\% confidence interval for $\pi(h)$ equal to:
\[I \ \equiv \ I_{100,000} \ = \ (0.0098, 0.0138)\]
This is again quite a small interval, of width $0.004$,
again providing good confidence about $\pi(h)$,
without needing to establish any CLT or difficult ergodicity property.

\ifx34
\section{Proof for Stationary Distribution}
Here we consider a more general case. Consider the Markov chain with
state space $\mathcal{X} = \{0, 1, 2, 3, \dots\}$, and transition
probability given by $p_{0, 0} = 1 - \delta$, $p_{0, 1} = \delta$,
and for all $n \geq 1$, $p_{n, n + 1} = \left(\frac{y}{y + 1}\right)^m$
and $p_{n, 0} = 1 - \left(\frac{y}{y + 1}\right)^m$. By the definition
of the stationary distribution, $$\pi(0) = (1 - \delta)\pi(0) + \sum_{i
= 1}^{\infty}(  1 - (\frac{i}{i + 1})^m) \pi(i)$$ $$\pi(1) = \delta
\pi(0) $$ $$\pi(x) =  (\frac{x-1}{x})^m\pi(x - 1), \forall x \geq 2.$$
The last two equations imply that $$\pi(x) = \frac{1}{x^m}\pi(1) =
\frac{\delta}{x^m}\pi(0), \ x \geq 1.$$  Check the first equation:
$$\sum_{i = 1}^{\infty} (1 - (\frac{i}{i + 1})^m)\pi(i) = \delta
\sum_{i = 1}^\infty (1 - (\frac{i}{i + 1})^m) \frac{1}{i^m}\pi(0) =
\delta \sum_{i = 1}^\infty (\frac{1}{i^m} - \frac{1}{(i + 1)^m}) \pi(0) =
\delta \pi(0).$$ So the first equation is satisfied. Finally we normalize
the measure: $$\sum_{i = 0}^{\infty}\pi(i) =  \pi(0) + \pi(0)\sum_{i =
1}^{\infty} \frac{\delta}{i^m} = \pi(0)(1 + \delta \zeta(m)),$$ where
\(\zeta(\cdot)\) is the Riemann zeta function. So $$\pi(0) = \frac{1}{1
+ \delta \zeta(m)}, \ \pi(x) = \frac{\delta}{x^m} \pi(0).$$ Finally,
substitute $\delta = 0.01$ and $m = 2$ into the formulas, we obtain
equation (1).
\fi

\subsection{A Polynomial-Tailed Non-Reversible Chain}

Define a Markov chain on the state space \(\mathcal{X} =
\mathbb{N}\), as follows.
Fix \(s > 2\). For all \(x > 1\), let \(p_{1, x} = x^{-s}\),
and \(p_{x, x-1} = 1\), with \(p_{x,y} = 0\) otherwise. Let \(\zeta(s)
= \sum_{i=1}^\infty i^{-s}\) denote the Riemann zeta function of
\(s\). It follows by induction that its stationary distribution $\pi$
satisfies that
\[\pi(x) = \left(\sum_{i=x}^\infty i^{-s}\right) \, \pi(1) \, , \quad x > 1.\]
Then, normalising the measure, we conclude that
\begin{align*}
\pi(1) &\ = \ \frac{1}{\zeta(s-1) - \zeta(s) + 1}\\
    \pi(x) &\ = \ \frac{\sum_{i=x}^\infty i^{-s}}{\zeta(s) - \zeta(s-1) +1}
\, , \quad x > 1.
\end{align*}

This chain is again non-reversible, and is again irreducible and aperiodic so
\(e_n\) will converge to \(\pi(h)\) whenever \(\pi(|h|) < \infty\).
The polynomial tails of $p_{x,1}$ and $\pi$ make
it difficult to directly establish a CLT, so we again
proceed through simulation.  Let $X_0=1$,
\(s = 5\), and \(h(x) = x\). For different numbers of
iterations \(n\), we run \(M = 1000\) replications of the chain
to estimate \(n\,\)Var\((e_n)\) and \(e_n\).

\begin{center}
\begin{tabular}{ |c|c|c|c| } 
\hline
\(n\) & \(n \, $\Big.$ \widehat{Var}(e_n)\) & \(\widehat{e_n}\) \\
\hline
100 & 0.552 & 1.61442 \\ 
200 & 0.807 & 1.61220 \\ 
500 & 1.013 & 1.61519 \\ 
1,000 & 0.946 & 1.61516 \\ 
2,000 & 0.874 & 1.61484 \\ 
5,000 & 0.883 & 1.61500 \\ 
10,000 & 0.774 & 1.61444\\ 
20,000 & 0.855 & 1.61436 \\ 
50,000 & 0.872 & 1.61441 \\ 
100,000 & 0.892 & 1.61458 \\ 
\hline
\end{tabular}
\end{center}
The simulations indicate that \(n\,\Var(e_n) \leq 1.5\) for all
\(n\). Therefore, we can take \(B = \sqrt{1.5}\). Then, setting \(\alpha
= 0.05\) and \(\epsilon = 0.001\), Theorem 1 gives an approximate conservative
\(95\%\) confidence interval for \(\pi(h)\) equal to:
\[I \equiv I_{100,000} = (1.597, 1.632)\]
which is again quite a small interval,
providing good confidence about $\pi(h)$.

\ifx34
\section{Derivation of the Stationary Distribution}
\noindent Based on the construction of the Markov chain, we know that the stationary distribution must satisfy the relation:
\begin{align*}
    \pi(1) &= (1 - \sum_{i=2}^\infty i^{-s}) \pi(1) + \pi(2) \\
    \pi(x) &= (x^{-s})\pi(1) + \pi(x+1)\\
\end{align*}

Using induction, it can be confirmed that
\[\pi(x) = (\sum_{i=x}^\infty i^{-s})\pi(1), \quad x > 1.\]
Lastly, we normalize the measure
\begin{align*}
    1 &= \sum_{i=1}^\infty \pi(i) \\
    &= \pi(1)(1 + \sum_{i=2}^\infty i^{-s} + \sum_{i=3}^\infty i^{-s} + \dots)\\
    &= \pi(1)(1 + 2^{-s} + 2(3^{-s}) + 3(4^{-s}) + \dots)\\
    &= \pi(1) (1 + \sum_{i=2}^\infty (i-1)(i^{-s}))\\
    &= \pi(1) (\sum_{i=1}^\infty i^{-(s-1)} - (\sum_{i=1}^\infty i^{-s} - 1))\\
\pi(1)    &= \frac{1}{\zeta(s-1) - \zeta(s) + 1}\\
    \pi(x) &= \frac{\sum_{i=x}^\infty i^{-s}}{\zeta(s) - \zeta(s-1) +1} \quad , x > 1.
\end{align*}
\fi

\subsection{A Discretised Irreversible Langevin Diffusion}


Finally, we examine a discrete version of the Irreversible Langevin
sampler for \(X_t \in \mathbb{R}^2\)
introduced in \cite{IrrevLang}.
Recall first the standard (reversible) Langevin diffusion, defined by:
\[dX_t \ = \ -\nabla U(X_t) \, dt + \sqrt{2D} \, dW_t,\]
where \(U : \mathbb{R}^2 \rightarrow \mathbb{R}\) is a $C^1$ function,
\(D > 0\) is a constant, and \(W_t\) is
two-dimensional Brownian motion. This process
converges to the stationary distribution having density
\[
\pi(x,y)
\ = \ \frac{e^{-U(x, y)/D}}{\int_{\mathbb{R}^2}e^{-U(x, y)/D} \, dx \, dy}.
\]
\cite{IrrevLang} show that the non-reversible family of diffusions
\[dX_t = \Big[-\nabla U(X_t) + C(X_t)\Big] \, dt + \sqrt{2D} \, dW_t,\]
also converges to this same stationary distribution $\pi$ provided that
the $C^1$ vector field \(C(x, y)\) satisfies
$\div(C e^{-2U}) =0$, and this condition is guaranteed if
\(C(x,y) = J \, \grad U(x, y)\) for some antisymmetric matrix $J$.
This new process is no longer reversible, but they argue that
it will sometimes converge faster.

We now consider
using this non-reversible Langevin diffusion
to estimate \(\pi(h)\) for some function
\(h:\mathbb{R}^2 \rightarrow \mathbb{R}\).
We run a discrete-time version of this continuous-time Irreversible
Langevin process, replacing each \(dt\) by $0.001$, and each
$dW_t$ by an independent $N(0,dt)$ draw.
We take
\(X_0 = (0, 0)\),
\(U(x, y) = (x^2 - 2)^2 + \frac{1}{4}y^2\),
$D=0.1$, 
\(C(x, y) = \begin{bmatrix}
0 & 10 \\
-10 & 0
\end{bmatrix} \, \grad U(x,y)\),
and \(h(x, y) = 16x^2 + 9y^2\).
We start each run at \(t=0\) and run until \(t=T\) for \(T \in \{10, 20,
\dots, 100\}\). Due to the complexity of the model,
we discard the results for \(0 \leq t
\leq 5\) as burn-in to minimize the effect of bias. Lastly, we define
\(n := (T-5)/dt\), and run \(M = 100\) replications for each \(T\) to
estimate \(e_n\) and \(n \, Var(e_n)\).
Our results are as follows:

\begin{center}
\begin{tabular}{ |c|c|c|c|c| } 
\hline
\(T\) & \(n\) & \(n \, $\Big.$ \widehat{Var}(e_n)\) & \(\widehat{e_n}\) \\
\hline
10 & 5,000 & 0.04666 & 32.00178 \\ 
20 & 15,000 & 0.02192 & 32.00156 \\ 
30 & 25,000 & 0.02074 & 32.00165 \\ 
40 & 35,000 & 0.02345 & 32.00173 \\ 
50 & 45,000 & 0.01570 & 32.00165 \\ 
60 & 55,000 & 0.01450 & 32.00164 \\
70 & 65,000 & 0.01668 & 32.00154 \\ 
80 & 75,000 & 0.01769 & 32.00168 \\ 
90 & 85,000 & 0.01612 & 32.00158 \\ 
100 & 95,000 & 0.01721 & 32.00172 \\ 
\hline
\end{tabular}
\end{center}
These simulations indicate that \(n\,\Var(e_n) \leq 0.05\)
for all \(n\). Therefore, we can take \(B = \sqrt{0.05}\). Then,
setting \(\alpha = 0.05\) and \(\epsilon = 0.001\), Theorem 1 gives
an approximate
conservative \(95\%\) confidence interval
for \(\pi(h)\) equal to:
\[I \equiv I_{95,000} =
(31.998, 32.005),\]
an extremely narrow interval (width $0.007$)
which gives good confidence about
the value of $\pi(h)$, again without proving any CLT or any challenging
ergodicity property.

\section{Summary}

In this paper, we have derived explicit asymptotic confidence
intervals for any MCMC algorithm with finite asymptotic variance,
started at any initial state, without requiring a Central Limit Theorem
nor reversibility nor any bias bound.  We have also derived explicit
non-asymptotic confidence intervals assuming bounds on the bias or first
moment, or alternatively that the chain starts in stationarity.  We have
related those non-asymptotic bounds to properties of MCMC bias, and
shown that polynomially ergodicity implies appropriate bias bounds.
Finally, we have applied our results to several numerical examples.
It is our hope that these results will provide simple and useful tools
for estimating errors of MCMC algorithms when CLTs are not easily available.

\bigskip\bigskip\noindent\bf Acknowledgements. \rm
We thank the anonymous referee for a very helpful report which led to
significant improvements in the manuscript.
JSR was partially supported by NSERC discovery grant RGPIN-2019-04142.
 
\bigskip


\bibliographystyle{apalike}
\bibliography{confbias}

\begin{thebibliography}{}

\bibitem[Alexopoulos and Seila, 1996]{10.1145/256562.256608}
Alexopoulos, C. and Seila, A.~F. (1996).
\newblock Implementing the batch means method in simulation experiments.
\newblock In {\em Proceedings of the 28th Conference on Winter Simulation}, WSC
  '96, page 214–221, USA. IEEE Computer Society.

\bibitem[Betancourt, 2020]{A1A2tweet}
Betancourt, M. (2020).
\newblock Reply tweet of {M}ay 25, 2020.
\newblock {\em Twitter}.
\newblock Available at:
  https://twitter.com/betanalpha/status/1264960932255600647.

\bibitem[Brooks et~al., 2011]{handbook}
Brooks, S., Gelman, A., Jones, G., and Meng, X.~L. (2011).
\newblock {\em {Handbook of Markov chain Monte Carlo}}.
\newblock CRC press.

\bibitem[Crane and Iglehart, 1975]{crane1975simulating}
Crane, M.~A. and Iglehart, D.~L. (1975).
\newblock Simulating stable stochastic systems: Iii. regenerative processes and
  discrete-event simulations.
\newblock {\em Operations Research}, 23(1):33--45.

\bibitem[Fishman, 1971]{fishman1971estimating}
Fishman, G.~S. (1971).
\newblock Estimating sample size in computing simulation experiments.
\newblock {\em Management Science}, 18(1):21--38.

\bibitem[Flegal and Jones, 2010]{FlegalJones}
Flegal, J.~M. and Jones, G.~L. (2010).
\newblock Batch means and spectral variance estimators in {M}arkov chain
  {M}onte {C}arlo.
\newblock {\em The Annals of Statistics}, 38(2):1034--1070.

\bibitem[Geyer, 1992]{geyer1992}
Geyer, C. (1992).
\newblock {Practical Markov chain Monte Carlo}.
\newblock {\em {Statistical Science}}, 7:473--483.

\bibitem[Geyer, 2011]{geyer2011introduction}
Geyer, C. (2011).
\newblock {Introduction to Markov chain Monte Carlo}.
\newblock {\em {Handbook of Markov chain Monte Carlo}}, 20116022:45.

\bibitem[Glynn and Whitt, 1972]{GlynnWhitt}
Glynn, P.~W. and Whitt, W. (1972).
\newblock The asymptotic validity of sequential stopping rules for stochastic
  simulations.
\newblock {\em The Annals of Applied Probability}, 2(1):180--198.

\bibitem[H{\"a}ggstr{\"o}m, 2005]{Haggstrom2005}
H{\"a}ggstr{\"o}m, O. (2005).
\newblock On the central limit theorem for geometrically ergodic {M}arkov
  chains.
\newblock {\em Probability Theory and Related Fields}, 132(1):74--82.

\bibitem[H{\"a}ggstr{\"o}m and Rosenthal, 2007]{haggstrom2007variance}
H{\"a}ggstr{\"o}m, O. and Rosenthal, J.~S. (2007).
\newblock {On variance conditions for Markov chain CLTs}.
\newblock {\em Electronic Communications in Probability}, 12:454--464.

\bibitem[Heidelberger and Welch, 1981]{10.1145/358598.358630}
Heidelberger, P. and Welch, P.~D. (1981).
\newblock A spectral method for confidence interval generation and run length
  control in simulations.
\newblock {\em Commun. ACM}, 24(4):233–245.

\bibitem[Hobert et~al., 2002]{HJPR}
Hobert, J., Jones, G., Presnell, B., and Rosenthal, J. (2002).
\newblock On the applicability of regenerative simulation in {M}arkov chain
  {M}onte {C}arlo.
\newblock {\em Biometrika}, 89:731--743.

\bibitem[Ibragimov and Linnik, 1971]{ibragimov1971independent}
Ibragimov, I.~A. and Linnik, Y.~V. (1971).
\newblock Independent and stationary sequences of random variables,(1971).
\newblock {\em Walters-Noordhoff, Gr{\"o}ningen}.

\bibitem[Jarner and Tweedie, 2003]{jarner2003}
Jarner, S.~F. and Tweedie, R.~L. (2003).
\newblock Necessary conditions for geometric and polynomial ergodicity of
  random-walk-type.
\newblock {\em Bernoulli}, 9(4):559--578.

\bibitem[Jones, 2004]{jones2004markov}
Jones, G.~L. (2004).
\newblock {On the Markov chain central limit theorem}.
\newblock {\em Probability Surveys}, 1:299--320.

\bibitem[Jones et~al., 2006]{JonesHaran}
Jones, G.~L., Haran, M., Caffo, B., and Neath, R. (2006).
\newblock Fixed-width output analysis for {M}arkov chain {M}onte {C}arlo.
\newblock {\em Journal of the American Statistical Association},
  101(476):1537--1547.

\bibitem[Jones and Hobert, 2001]{jones2001honest}
Jones, G.~L. and Hobert, J.~P. (2001).
\newblock {Honest exploration of intractable probability distributions via
  Markov chain Monte Carlo}.
\newblock {\em Statistical Science}, pages 312--334.

\bibitem[Kipnis and Varadhan, 1986]{kipnis1986central}
Kipnis, C. and Varadhan, S. R.~S. (1986).
\newblock Central limit theorem for additive functionals of reversible {M}arkov
  processes and applications to simple exclusions.
\newblock {\em Communications in Mathematical Physics}, 104(1):1--19.

\bibitem[{L}atuszy{\'n}ski et~al., 2013]{latuszynski2013nonasymptotic}
{L}atuszy{\'n}ski, K., Miasojedow, B., and Niemiro, W. (2013).
\newblock Nonasymptotic bounds on the estimation error of {MCMC} algorithms.
\newblock {\em Bernoulli}, 19(5A):2033--2066.

\bibitem[Law and Kelton, 1984]{10.2307/170945}
Law, A.~M. and Kelton, W.~D. (1984).
\newblock Confidence intervals for steady-state simulations: I. a survey of
  fixed sample size procedures.
\newblock {\em Operations Research}, 32(6):1221--1239.

\bibitem[Mengersen et~al., 1999]{mengersen1999mcmc}
Mengersen, K.~L., Robert, C.~P., and Guihenneuc-Jouyaux, C. (1999).
\newblock {MCMC} convergence diagnostics: a review.
\newblock {\em Bayesian Statistics}, 6:415--440.

\bibitem[Rey-Bellet and Spiliopoulos, 2015]{IrrevLang}
Rey-Bellet, L. and Spiliopoulos, K. (2015).
\newblock Irreversible {L}angevin samplers and variance reduction: a large
  deviations approach.
\newblock {\em Nonlinearity (Print)}, 28(7).

\bibitem[Roberts and Rosenthal, 2004]{roberts2004general}
Roberts, G.~O. and Rosenthal, J.~S. (2004).
\newblock {General state space Markov chains and MCMC algorithms}.
\newblock {\em Probability Surveys}, 1:20--71.

\bibitem[Rosenthal, 2017]{rosenthal2017}
Rosenthal, J.~S. (2017).
\newblock {Simple confidence intervals for MCMC without CLTs}.
\newblock {\em Electron. J. Statist.}, 11(1):211--214.

\bibitem[T{\'o}th, 1986]{Toth1986}
T{\'o}th, B. (1986).
\newblock Persistent random walks in random environment.
\newblock {\em Probability Theory and Related Fields}, 71(4):615--625.

\end{thebibliography}

\end{document}